\documentclass[10pt,twoside]{article}

\usepackage{graphicx}
\usepackage [usenames] {color}
\usepackage [colorlinks=true,
linkcolor=webgreen,
filecolor=webbrown,
citecolor=webgreen] {hyperref}
\definecolor {webgreen} {rgb} {0,.5,0}
\definecolor {webbrown} {rgb} {.6,0,0}
      
\usepackage {color}
\usepackage {float}

\usepackage{amsfonts}
\usepackage{amstext}
\usepackage{amssymb,amsmath,amsbsy}

\usepackage{hyperref}
\newtheorem{theorem}{Theorem}

\newtheorem{definition}{Definition}
\newtheorem{lemma}{Lemma}
\newtheorem{proposition}{Proposition}

\newcounter{thcount}
\def\R{\mathbb{R}}

\def\E{\mathbb{E}}
\def\D{\mathbb{D}}
\def\WW{\mathbb{W}}
\def\0{\mathbf{0}}
\def\1{\mathbf{1}}

\def\d{\mathbf{d}}

\def\x{\mathbf{x}}

\def\u{\mathbf{u}}

\def\y{\mathbf{y}}

\def\I{\mathbf{I}}

\def\A{\mathbf{A}}
\def\M{\mathbf{M}}
\def\L{\mathbf{L}}
\def\W{\mathbf{W}}
\def\DD{\mathbf{D}}
\def\NN{\mathbf{N}}
\def\PP{\mathbf{P}}

\def\tr{\mathrm{tr}}
\def\diag{\mathrm{diag}}

\def\rk{\mathrm{rank}}
\def\Var{{\mathrm{Var}}}
\def\Vol{{\mathrm{Vol}}}
\def\cov{{\mathrm{Cov}}}
\def\corr{{\mathrm{Corr}}}

\def\det{{\mathrm{det}}}

\def\part{\cal P} 

\begin{document}

%{\Large {\it Research Article}}\\

{\Large {\bf {When the largest eigenvalue of the modularity and the normalized
modularity matrix is zero}}\\

{\large {\bf Marianna Bolla, Brian Bullins, Sorathan Chaturapruek, Shiwen Chen,
Katalin Friedl}}\\

{Institute of Mathematics, Budapest University of Technology and Economics 
and Budapest Semester of Mathematics}

%Correspondence should be addressed to Marianna~Bolla, 
%{\tt {marib@math.bme.hu}}\\

\date{}

\renewcommand\abstractname{Abstract}

\begin{abstract}
\noindent \\
In July 2012, at the Conference on Applications of Graph Spectra in Computer 
Science, Barcelona, D. Stevanovic posed the following open problem:
which graphs have the zero as the largest eigenvalue of their modularity matrix?
The conjecture was that only the complete and the complete multipartite graphs.
They indeed have this property, but are they the only ones? 
In this paper, we will give
an affirmative answer to this question and prove a bit more: both the modularity
and the normalized modularity matrix of a graph is negative semidefinite if and
only if the graph is complete or complete multipartite.


\end{abstract}

\section {Introduction}
\label{intro}

In \cite{New1} Newman and Girvan defined the modularity matrix of a simple graph
on $n$ vertices with an $n\times n$ symmetric adjacency matrix $\A$ as
$$
 \M =\A -\frac1{2e} \d \d^T ,
$$
where $\d =(d_1 ,\dots ,d_n )^T$ is the so-called \textit{degree-vector}
comprised of the vertex-degrees $d_i$'s and $2e =\sum_{i=1}^n d_i$ is twice
the number of edges. In \cite{Bol6} we formulated the modularity matrix of an
\textit{edge-weighted graph}  $G=(V, \W)$ on the
$n$-element vertex-set $V$ with an $n\times n$ symmetric weight-matrix 
$\W$ -- the  entries of which are pairwise similarities 
between the vertices and 
satisfy $w_{ij}=w_{ji}\ge 0$, $w_{ii}=0$ 
$(i=1,\dots ,n )$, further, $\sum_{i=1}^n \sum_{j=1}^n w_{ij}=1$ -- as follows.
$$
 \M =\W -\d \d^T ,
$$
where the entries of $\d$ are now the \textit{generalized vertex-degrees}
$d_i =\sum_{j=1}^n w_{ij}$ $(i=1,\dots ,n )$. The assumption 
$\sum_{i=1}^n d_{i} =1$ does not hurt the generality, but simplifies 
further notation and makes it possible to consider $\W$ as a symmetric
joint distribution of two identically distributed discrete random variables
taking on $n$ different values.

The modularity matrix $\M$ always has a zero eigenvalue
with eigenvector $\1 =\1_n=(1,\dots ,1)^T$, since its rows sum to zero.
Because $\tr (\M )<0$, $\M$ must have at least one negative eigenvalue, and
it is usually
indefinite. For the complete and the complete multipartite graphs, however, its
largest eigenvalue is zero, as we will show in Section~\ref{ex}. 
In Theorem~\ref{th1} of Section~\ref{th}, 
we will prove that the modularity matrix of a 
simple graph is negative semidefinite if and
only if it is complete or complete multipartite. In Theorem~\ref{th2}
we will extend
this statement to the negative semidefiniteness of the normalized modularity
matrix introduced in \cite{Bol6} as 
$$
  \M_D =\DD^{-1/2} \M \DD^{-1/2} , 
$$
where $\DD =\diag (d_1 ,\dots ,d_n )$ is the \textit{degree-matrix}. 
The eigenvalues of $\M_D$ are the same, irrespective  of
whether we start with the adjacency
or normalized edge-weight matrix of a simple graph, and they are in the
$[-1, 1]$ interval. $\M_D$ is closely related to the normalized Laplacian;
therefore, our statements have important consequences, as
 established in Section~\ref{pre},
 for the normalized Laplacian spectrum.

In Section~\ref{fin}, we discuss some other implications of 
Theorems~\ref{th1} and~\ref{th2} 
concerning the Newman--Girvan modularity of \cite{New1,New5} and the
maximal correlation of \cite{Bolla1}.

\section{Preliminaries}
\label{pre}

First we introduce some notation.
\begin{definition}
The simple graph on $n$ vertices is complete if the entries of its adjacency
matrix are
$$ a_{ij} := \left\{
      \begin{array}{ll}   1 & \mbox{if } \quad i\ne j \\
     0   & \mbox{if } \quad i=j .
    \end{array}
\right.
$$ 
This graph is denoted by $K_n$.
\end{definition}

\begin{definition}
The simple graph on the  $n$-element vertex-set $V$ is complete multipartite 
with $2\le k \le n$ 
partites (color-classes) $V_1 ,\dots ,V_k$ (they form a partition of the vertices) 
if the entries of its adjacency matrix are
$$ a_{ij} := \left\{
      \begin{array}{ll}   1 & \mbox{if } \quad c(i)\ne c(j) \\
     0   & \mbox{if } \quad c(i)=c(j) , 
    \end{array}
\right.
$$ 
where $c(i)$ is the color of vertex $i$. Here the non-empty,
disjoint vertex-subsets form so-called maximal independent sets of the 
vertices. If $|V_i |=n_i$ $(i=1,\dots k)$, $\sum_{i=1}^k n_i =n$, then  
this graph is denoted by $K_{n_1 ,\dots ,n_k}$.
\end{definition}

Note that $K_n$ is also complete multipartite with $n$ partites, i.e., 
it is the $K_{1,\dots ,1}$ graph; therefore, in the sequel, whenever we speak of
complete multipartite graphs, complete graphs are also understood.

In \cite{Bolla} we introduced the \textit{normalized Laplacian} of $G=(V, \W)$ 
as $\L_D =\I -\DD^{-1/2} \W \DD^{-1/2}$, and in \cite{Bol6,Bol7} we established
the following relation between the spectra of $\L_D$ and $\M_D$ when
$G$ is connected ($\W$ is irreducible).
Let $0=\lambda_0 <\lambda_1 \le \dots \le \lambda_{n-1} \le 2$  denote
eigenvalues of $\L_D$ with 
corresponding
unit-norm, pairwise orthogonal eigenvectors $\u_0 ,\dots ,\u_{n-1}$.
Namely, $\u_0 =(\sqrt{d_1},\dots , \sqrt{d_n})^T$, which will be
denoted by  $\sqrt{\d }$. 
Enumerating  the eigenvalues of $\M_D $ of the same connected graph in the 
order 
%\I -\L_D -\sqrt{\d} \sqrt{\d}^T$
$1>\mu_1 \ge \dots \ge \mu_{n-1} \ge -1$, we have
$\mu_i = 1-\lambda_i$ with the same 
eigenvector $\u_i$  $(i=1,\dots ,n-1)$; further, $\mu_n =0$ with 
corresponding unit-norm eigenvector $\sqrt{\d}$.

The smallest positive normalized Laplacian eigenvalue $\lambda_1$ solves the
following quadratic placement problem:
\begin{equation}\label{mini}
 \lambda_1 = \min \sum_{i<j} w_{ij} (r_i -r_j )^2 
\end{equation}
subject to
\begin{equation}\label{const}
 \sum_{i=1}^n d_i  r_i =0 \quad \textrm{and} \quad\sum_{i=1}^n d_i  r_i^2 =1 .
\end{equation}
In \cite{Bolla} we proved that
the optimal vertex-representatives  $r_1 ,\dots ,r_n$ giving the above 
minimum are the coordinates of the vector $\DD^{-1/2} \u_1$.

\begin{proposition}\label{prop} 
Let $G =(V, \W )$ be an edge-weighted graph, the 
weight-matrix of which has at least one off-diagonal zero entry. Then 
$\lambda_1 (G) \le 1$.
\end{proposition}

\textbf{Proof}:
Since the normalized Laplacian spectrum of isomorphic graphs is the same
(irrespective of the same permutation of the rows and columns of $\W$), we
may assume that $w_{12} =0$. Let us define the following representation of
the vertices:
$$
 r_1 := \frac{d_2}{\sqrt{d_2^2 d_1 +d_1^2 d_2}}, \quad
 r_2 :=-\frac{d_1}{\sqrt{d_2^2 d_1 +d_1^2 d_2}} ,
$$
and $r_i :=0$, $i=3,\dots ,n$ when $n$ exceeds 2. These representatives
satisfy conditions (\ref{const}) and in view of (\ref{mini}):
$$
 \lambda_1 (G)\le \sum_{i<j} (r_i -r_j )^2 w_{ij} =
 \frac{ \sum_{j\ne 1} (d_2 -0)^2 w_{1j} +\sum_{j\ne 2} (-d_1 -0)^2 w_{2j} }
  {d_2^2 d_1 +d_1^2 d_2} =1 ,
 $$
which finishes the proof.

We know (see e.g., \cite{Bol13}) that %for the complete graph,
$\lambda_1 (K_n )=\dots =\lambda_{n-1} (K_n)=\frac{n}{n-1}$. 
Proposition~\ref{prop} guarantees that all the other simple graphs have
$\lambda_1 \le 1$. This was also proved in \cite{Chung0}.
In Section~\ref{th} we will 
prove that equality is attained only for $K_{n_1 ,\dots ,n_k}$ $(k<n)$.

\section{Modularity spectra of complete and complete multipartite 
graphs}
\label{ex}

Here we calculate modularity spectra of the exceptional graphs in question.
\begin{proposition}\label{ex1}
The spectrum of $\M (K_n )$ consists of the single
eigenvalue 0 with eigenvector $\1_n$ and the number $-1$ with multiplicity
$n-1$ and eigen-subspace $\1_n^{\perp }$.
\end{proposition}

\textbf{Proof}:
The adjacency matrix of $K_n$ is 
$\A (K_n )=\1_n \1_n^T -\I_n $, $\d =(n-1) \1_n$, $2e =n(n-1)$, hence
$$
\begin{aligned}
 \M (K_n ) &= \frac1{n} \1_n \1_n^T -\I_n =
 (\frac{\1}{\sqrt{n}}) (\frac{\1}{\sqrt{n}})^T -
 \left[ (\frac{\1}{\sqrt{n}}) (\frac{\1}{\sqrt{n}})^T +
   \sum_{i=2}^n 1\cdot\u_i \u_i^T \right]  \\
 &=\sum_{i=2}^n (-1)\cdot\u_i \u_i^T ,
\end{aligned}
$$
where $\u_2 ,\dots \u_n$ is an arbitrary orthonormal set in $\1_n^{\perp}$.
Therefore, the unique spectral decomposition of $\M (K_n )$ is as stated
in the proposition.

\begin{proposition}\label{ex2}
The spectrum of $\M_D (K_n )$ consists of the single
eigenvalue 0 with eigenvector $\sqrt{\d }$ and the number $-\frac{1}{n-1}$ 
with multiplicity $n-1$ and eigen-subspace $\sqrt{\d}^{\perp }$.
\end{proposition}

This proposition follows from the characterization of the normalized
Laplacian spectrum of $K_n$ given in \cite{Bol13}.

\begin{proposition}\label{ex3}
The spectrum of $\M (K_{n_1 ,\dots ,n_k})$ consists of $k-1$ strictly negative
eigenvalues and zero with multiplicity $n-k+1$. 
\end{proposition}

\textbf{Proof}:
The adjacency matrix $\M$ 
of the complete multipartite graph $K_{n_1 ,\dots ,n_k}$ is
a block-matrix with diagonal blocks of all zeros and off-diagonal 
blocks of all 1's. Let $V_1 ,\dots ,V_k$ denote the independent, disjoint 
vertex-subsets, $|V_i | =n_i$, $i=1,\dots ,k$;
$d_l = n-n_i$ if $l\in V_i$; $2e =\sum_{l=1}^n d_l =\sum_{i=1}^k n_i (n-n_i )
=n^2 - \sum_{i=1}^k n_i^2$. Therefore, $\M$ is a blown-up matrix (see 
\cite{Bol13}) with blow-up 
sizes $n_1 ,\dots ,n_k$ of the $k\times k$ pattern matrix $\PP$ with entries
$$
 p_{ij} = (1-\delta_{ij}) -\frac{(n-n_i)(n-n_j )}{2e} ,
$$
where $\delta_{ij}$ is the Kronecker-delta. 
Consequently, $\rk (\M ) =\rk (\PP ) \le k$. 
We will prove that its rank is $k-1$,
it has $k-1$ negative eigenvalues, and all its other eigenvalues are zeros. 
An eigenvector $\u$ belonging to a nonzero eigenvalue $\lambda$ is piecewise 
constant with $n_1$ coordinates equal to $y_1$,\dots ,
$n_k$ coordinates equal to $y_k$. With these, the eigenvalue--eigenvector
equation yields  that 
$$
 \sum_{j=1}^k n_j \left[ (1-\delta_{ij}) -\frac{(n-n_i)(n-n_j )}{2e} \right]
  y_j =\lambda y_i  .
$$ 
Therefore, $\lambda$ is an eigenvalue of the $k\times k$ matrix $\PP \NN$
with eigenvector $(y_1 ,\dots ,y_k )^T$, where $\NN =\diag (n_1 ,\dots ,n_k )$.
The matrix $\PP \NN$ is not symmetric, but its eigenvalues are real because of 
the
above, or else its eigenvalues are also eigenvalues of the symmetric matrix 
$\NN^{1/2} \PP \NN^{1/2}$. It is easy to see that the row sums of $\PP \NN$
are zeros:
$$
 \sum_{j=1}^k n_j \left[ (1-\delta_{ij}) -\frac{(n-n_i)(n-n_j )}{2e} \right]=0,
$$
i.e.,
\begin{equation}\label{eigeneq}
\sum_{j=1}^k n_j p_{ij}
  y_j =\lambda y_i  .
\end{equation}
Therefore, zero is an eigenvalue with eigenvector $\1_k$, which gives another 
zero eigenvalue of $\M$ with eigenvector $\1_n$. Thus, zero is an eigenvalue
of $\M$ with multiplicity $n-k+1$ and corresponding eigensubspace of this
dimension, including $\1_n$. 

Now we will prove that all the non-zero eigenvalues are negative. In view of
(\ref{eigeneq}),
$$
 \lambda \sum_{i=1}^k n_i y_i = \sum_{i=1}^k n_i (\lambda y_i ) =
 \sum_{j=1}^k n_j y_j \sum_{i=1}^k n_i p_{ij} =0 .
$$
Consequently, if $\lambda\ne 0$, then $\sum_{i=1}^k n_i y_i =0$.
On the other hand,
$$
 \lambda \sum_{i=1}^k n_i y_{i}^2 =\sum_{i=1}^k (n_i y_{i}) (\lambda y_i ) =
\sum_{i=1}^k n_i y_{i} \sum_{j=1}^k  p_{ij} n_j y_j =
\sum_{i=1}^k \sum_{j=1}^k p_{ij} (n_i y_{i}) (n_j y_j ) .
$$
We will show that the right hand side is negative, and therefore, by
$\sum_{i=1}^k n_i y_{i}^2 >0$, we get that $\lambda <0$. Indeed,
$$
\begin{aligned}
&\sum_{i=1}^k \sum_{j=1}^k p_{ij} (n_i y_{i}) (n_j y_j ) =
\sum_{i=1}^k \sum_{j=1}^k  \left[ (1-\delta_{ij}) -\frac{(n-n_i)(n-n_j )}{2e}
 \right] (n_i y_{i}) (n_j y_j )  \\
 &=  \sum_{i=1}^k \sum_{j=1}^k  (1-\delta_{ij})  (n_i y_{i}) (n_j y_j ) -
 \frac1{2e} \left[ \sum_{i=1}^k (n-n_i) n_i y_i \right]
   \left[ \sum_{j=1}^k (n-n_j) n_j y_j \right] \\
 &= \left( \sum_{i=1}^k n_i y_i \right) \left( \sum_{j=1}^k n_j y_j \right) -
  \sum_{i=1}^k (n_i y_i )^2 -
   \frac1{2e} \left[ \sum_{i=1}^k (n-n_i )n_i y_i \right]^2 <0 ,
\end{aligned}
$$
which, by $\sum_{i=1}^k n_i y_i =0$, finishes the proof.

\begin{proposition}\label{ex4}
$\M_D (K_{n_1 ,\dots ,n_k})$ is also negative semidefinite.
\end{proposition}

\textbf{Proof}:
We have to show that for any $\x \in \R^n$, $\x^T \M_D  \x$
is nonpositive, where for brevity, $\M_D$ denotes the normalized
modularity matrix of  $K_{n_1 ,\dots ,n_k}$. In fact,
$$
 \x^T \M_D \x = (\DD^{-1/2}\x )^T \M (\DD^{-1/2} \x ) =\y^T \M \y \le 0
$$
for any $\y\in \R^n$ because of the negative semidefiniteness of the
modularity matrix of $K_{n_1 ,\dots ,n_k}$. Due to the invertibility of
the degree-matrix $\DD$ (our graph is connected, hence cannot have isolated 
vertices), the above relation holds for any $\x$ as well.

Furthermore, $\M_D (K_{n_1 ,\dots ,n_k})$ has rank 
$k-1$ with $k-1$ strictly negative eigenvalues, and the  
$(n-k+1)$-dimensional eigensubspace corresponding to the zero
eigenvalue looks like
\begin{equation}\label{alter}
 \{ \x : \, \sum_{j \in V_i} \sqrt{d_j} x_j =0, \, i=1,\dots ,k \} =
  \{ \x : \, \sum_{j \in V_i} x_j =0, \, i=1,\dots ,k \} .
\end{equation}
We can use that the vertex-degrees within the partites are the same,
and so, we have a blown-up matrix again.

Note that in the case of $k=n$, the results of Propositions~\ref{ex3} and 
\ref{ex4} exactly provide those of Propositions~\ref{ex1} and 
\ref{ex2}.

\section{The main results and proofs}
\label{th}

To prove the main results, we will intensively use the following 
characterization of complete multipartite graphs, including the complete graph.
Although this is a known result of graph theory, we enclose the proof as well.

\begin{lemma}\label{lem}
A simple connected graph is complete multipartite if and only if it
has no 3-vertex induced subgraph with exactly one edge. 
\end{lemma}

\textbf{Proof}:
We will call the above subgraph \textit{forbidden pattern}.

\begin{itemize}
\item In the forward direction, 
a complete multipartite graph can have the following 
types of 3-vertex induced subgraphs (not all of them appear necessarily,
only if the size of partites allows it):
\begin{itemize}
\item the three vertices are from the same partite, 
in which case the induced subgraph has no edges;
\item the three vertices are from three different partites, 
in which case the induced subgraph is the complete graph $K_3$;
\item two of the vertices are from the same, and the third from a different
partite, in which case the induced subgraph has exactly two edges 
(called cherry).
\end{itemize}
None of them is the forbidden pattern.

\item Conversely, suppose that our graph does not have the forbidden pattern.
The following procedure shows that it is then complete multipartite. 
Let the first
cluster be a maximal independent set of the vertices, say 
$V_1$. We claim that each vertex in ${\overline V}_1$ is 
connected to each vertex of $V_1$. Indeed, let $c\in {\overline V}_1$ be a 
vertex;
it must be connected to a vertex (say, $a$) of $V_1$, since if not, it 
could be joined to $V_1$, which contradicts the 
maximality of $V_1$ as an
independent set.  If $c$ were not connected to another $b\in V_1$, then 
%for some $a\in V_1$: $a\sim c$, and for an other $b\in V_1$: $b\nsim c$,
$a,b,c$ would form a forbidden pattern, but our graph does not contain such
in view of our starting assumption.

Then let $V_2$ be a  maximal independent set of vertices within 
${\overline V}_1$, 
say $V_2$.  We claim that each vertex in $\overline{V_1 \cup V_2}$ 
is connected to each vertex of $V_1$ and
$V_2$. The connectedness to vertices of $V_1$ is already settled. 
By the maximality of $V_2$ as an independent set, any vertex of
$\overline{V_1 \cup V_2}$ must be connected to at least one vertex of $V_2$.
If we found 
a vertex $c\in \overline{V_1 \cup V_2}$ such that
for some $a\in V_2$: $a\sim c$, and for another $b\in V_2$: $b\nsim c$,
then $a,b,c$ would form a forbidden pattern, which is excluded. 

Advancing in this way, one can see that
the procedure produces maximal disjoint independent sets of vertices such that
the independent vertices of $V_k$ are connected to every vertex in $V_1 ,\dots
,V_{k-1}$. At each step we can select a maximal independent set out of the
remaining vertices; in the worst case it contains only one vertex.
The absence of the forbidden 
pattern guarantees that we can always continue our algorithm until all vertices
are placed into a cluster.
This procedure will exhaust the set of vertices and result in a
complete multipartite graph. 
The point is that in the absence of the forbidden pattern we can divide the
vertices into independent sets which are fully connected.
\end{itemize}

If we proceed with non-increasing cardinalities of $V_i$'s, then one-vertex
independent sets may emerge at the end of the process. 
Moreover, up to the labeling of the vertices and
the numbering of the independent sets, 
the resulting multipartite structure is unique.  
In fact, the above procedure just recovers this unique structure
in the absence of the forbidden pattern.

Now, the answer to the open question follows.
\begin{theorem}\label{th1}
 The modularity matrix of a 
simple connected graph is negative semidefinite if and
only if it is complete multipartite. 
\end{theorem}

\textbf{Proof}:
First we prove that when a simple graph is not complete multipartite, then
its modularity matrix cannot be negative semidefinite.  
By Lemma~\ref{lem}, a simple graph is not complete multipartite if and only
if it contains the forbidden pattern. Let us take such a graph. 
Since the modularity spectrum does not
depend on the labeling of the vertices, assume that the first three vertices
form the forbidden pattern, i.e., the upper left corner of the adjacency matrix
is 
$$                                          
\begin{pmatrix}  0 & 1 & 0 \\
                 1 & 0 & 0 \\
                 0 & 0 & 0
      \end{pmatrix} .
$$
It is known that a matrix is negative semidefinite if and only if its every 
principal minor of odd order is non-positive, and every 
principal minor of even order is non-negative. Since the principal minor
of order 3 of this graph's  modularity matrix is
$$
 \left( \frac1{2e} \right)^3 \det \begin{pmatrix}  
                -d_1^2 & 2e -d_1 d_2 & -d_1 d_3 \\
                 2e -d_1 d_2 & -d_2^2 & -d_2 d_3 \\
                 -d_1 d_3  & -d_2 d_3  & -d_3^2
      \end{pmatrix}  = \frac1{8e^3} 4e^2 d_3^2 = \frac{d_3^2}{2e} >0 ,
$$
the modularity matrix cannot be negative semidefinite. This fact, together
with Proposition~\ref{ex3}, finishes the proof.

We are able to prove a similar statement for the normalized modularity matrix.
\begin{theorem}\label{th2}
The normalized modularity matrix of a 
simple connected graph is negative semidefinite if and
only if it is complete multipartite. 
\end{theorem}

\textbf{Proof}: Now we will prove that if a simple graph is not
complete multipartite, or equivalently, if it contains 
the forbidden pattern, then the largest eigenvalue  of its 
normalized modularity matrix is strictly  positive. This fact,
together with Proposition~\ref{ex4}, will finish the proof.

Referring to \cite{Bol7}, the largest eigenvalue $\mu_1$ of $\M_D$ is the
second largest eigenvalue of $\DD^{-1/2} \W \DD^{-1/2}$, whose largest eigenvalue
is 1 with corresponding eigenvector $\sqrt{\d }$ (this is unique if our
graph is connected). Therefore, we think in terms of the two largest eigenvalues
of  $\DD^{-1/2} \W \DD^{-1/2}$.
We can again assume that the first three vertices form the forbidden
pattern and so, the upper left corner of this matrix looks like
$$                                          
\begin{pmatrix}  0 & \frac{1}{\sqrt{d_1 d_2}} & 0 \\
                 \frac{1}{\sqrt{d_1 d_2}} & 0 & 0 \\
                 0 & 0 & 0
      \end{pmatrix} .
$$
 Then by the Courant--Fischer--Weyl minimax principle:
$$
\mu_1 =\max_{\substack{ \| \x \| =1 \\ \x^T \sqrt{\d} =0 }}
 \x^T \DD^{-1/2} \W \DD^{-1/2} \x .
$$
Therefore, to prove that $\mu_1 >0$, it suffices to find an $\x \in \R^n$ 
($n$ is the number of vertices) that 
satisfies conditions $\| \x \| =1$, $\x^T \sqrt{\d} =0$ and for which
$\x^T \DD^{-1/2} \W \DD^{-1/2} \x >0$. 
(The unit norm condition can be relaxed here, because $\x$ can later be
normalized, without changing the sign of the above quadratic form.)

Indeed, let us look for $\x$ in the form $\x = (x_1 ,x_2,x_3 ,0,\dots ,0)^T$
such that 
\begin{equation}\label{felt}
\sqrt{d_1 } x_1 +\sqrt{d_2 } x_2 +\sqrt{d_3 } x_3 =0 .
\end{equation}
Then the inequality
$$
 \x^T \M_D \x = \frac{2x_1 x_2}{\sqrt{d_1 d_2}} >0
$$
can be satisfied with any 
$\x =(x_1 ,x_2,x_3 ,0,\dots ,0)^T$ such
that $x_1$ and $x_2$ are both positive or both negative, and due to
(\ref{felt}),  
$$
 x_3 =-\frac{\sqrt{d_1 } x_1 +\sqrt{d_2 } x_2}{\sqrt{d_3 }}
$$
is a good choice, which will have the opposite sign.
(Note that $d_i$'s are positive, since we deal with connected graphs.)

\section{Conclusions}
\label{fin}

The results of Section~\ref{th} have the following important implications.
\begin{itemize}
\item
In terms of $\mu_1$, a result of \cite{Bolla1} can be interpreted in the
following way.  
We use the setup of correspondence analysis, applied to the symmetric joint
distribution embodied by the entries of $\W$.
Let $\psi$ and $\phi$ be identically distributed (i.d.) random variables
with this joint distribution.  
%(they are usually not independent). 
Say, these  discrete random variables take on values 
$r_1 ,\dots ,r_n$ with probabilities $d_1 ,\dots ,d_n$ (margin of the
joint distribution; the two margins are the same, since $\W$ is symmetric).
%and this defines a probability distribution, since the sum of the entries
%of $\W$ is 1).
Then
$$
 \mu_1 =\max_{\psi ,\phi \, \textrm{i.d.} } \corr_{\WW } (\psi , \phi ) =
 \max_{\substack{ \psi ,\phi \, \textrm{i.d.} \\
  \Var_{\D } \psi =1 \\ \E_{\D } \psi =0 }} 
  \cov_{\WW } (\psi , \phi ) =
 \max_{\substack{ \sum_{i=1}^n d_i  r_i =0 \\ \sum_{i=1}^n d_i  r_i^2 =1 }}
  \sum_{i=1}^n \sum_{i=1}^n w_{ij} r_i r_j ,
$$
and the maximum is attained when the values $r_1 ,\dots ,r_n$ are coordinates 
of the vector $\DD^{-1/2} \u_1$. 
In this setup, the conditions for the zero expectation and
unit variance are analogous to those of (\ref{const}).
In \cite{Bolla1} $\mu_1$ is called \textit{symmetric maximal correlation}. 
The results of the present paper show that it is positive if and only if the 
joint distribution is not of a complete multipartite structure.  

\item
The $2$-way Newman--Girvan 
modularity (see \cite{New1,New5,Bol6}) of $G=(V,\W)$ is
$$
 Q_2 = \max_{\emptyset \ne U\subset V} Q( U, {\overline U}) ,
$$
where the modularity of the $2$-partition $(U,{\overline U})$ of $V$ is written
in terms of the entries $m_{ij}$'s (summing to 0) of $\M (G)$:
$$
\begin{aligned}
 Q( U, {\overline U}) &= 
 \sum_{i,j\in U} m_{ij} + 
 \sum_{i,j\in {\overline U}} m_{ij} =-2\sum_{i\in U, \, j\in {\overline U}} m_{ij} \\
 &= - 2 [w (U ,{\overline U} ) -\Vol (U) \Vol ({\overline U}) ] , 
\end{aligned}
$$
where   $w (U ,{\overline U} )= \sum_{i\in U} \sum_{j\in {\overline U}} w_{ij}$
is the weighted cut between $U$ and $\overline U$, whereas
$\Vol (U) =\sum_{i\in U} d_i$ is the volume of the vertex-subset $U$.
These formulas are valid under the condition $\Vol (V)=1$; otherwise, they
should be adjusted by $2e$.

Now we use the idea of the proof of the 
Expander Mixing Lemma extended to edge-weighted graphs (see \cite{Bol7}). 

With the notation of Section~\ref{pre} and introducing $\mu_0 =1$, $\u_0
=\sqrt {\d}$, 
$$
 \DD^{-1/2} \W \DD^{-1/2} =\sum_{i=0}^{n-1}\mu_i \u_i \u_i^T
$$ 
is spectral decomposition.

Let $U\subset V$  be arbitrary and the indicator vector of $U$ is
denoted by $\1_U \in \R^n$. Further,
put $\x :=\DD^{1/2} \1_U$ and $\y :=\DD^{1/2}\1_{\overline U}$, and
let $\x = \sum_{i=0}^{n-1} a_i \u_i$ and
$\y = \sum_{i=0}^{n-1} b_i \u_i$ be the expansions of $\x$ and $\y$ in the 
orthonormal basis $\u_0 ,\dots ,\u_{n-1}$ with coordinates
$a_i =\x^T \u_i$ and $b_i =\y^T \u_i$, respectively. Observe that
$w (U,{\overline U})= \1_U^T \W \1_{\overline U} =\x^T (\DD^{-1/2} \W \DD^{-1/2}) \y^T$ and
$\1_{\overline U} = \1_n -\1_U$; therefore,
$$
 b_i =\y^T \u_i =\DD^{1/2} (\1 -\1_U )\u_i = \u_0^T \u_i -\x^T \u_i =-a_i  \,
(i=1,2,\dots ,n-1).
$$ 
Further, $a_0 =\Vol (U)$ and $b_0 =\Vol ({\overline U})$.
%$\sum_{i=0}^{n-1} a_i^2 =\| \x \|^2 =\Vol (X)$,
%$\sum_{i=0}^{n-1} b_i^2 =\| \y \|^2 =\Vol (Y)$.
Based on these observations,
$$
 w (U, {\overline U}) - \Vol (U) \Vol ({\overline U}) =
 \sum_{i=1}^{n-1} \mu_i a_i b_i  = -\sum_{i=1}^{n-1} \mu_i a_i^2 .
$$
Consequently,  
$Q(U, {\overline U}) =2\sum_{i=1}^{n-1} \mu_i a_i^2$. Therefore, provided 
that the normalized
modularity matrix of the underlying graph is negative semidefinite 
(or equivalently, our graph is complete multipartite),
$Q(U, {\overline U}) \le 0$
for all 2-partitions of the vertices, and hence, the 
 $2$-way Newman--Girvan 
modularity is also non-positive (in most cases, it is negative). 
Nonetheless this property
does not characterize the complete multipartite graphs. There are graphs
with positive $\mu_1$ and zero or sometimes negative 
 $2$-way Newman--Girvan modularity.

\item
Recall that the smallest positive normalized Laplacian eigenvalue 
$\lambda_1$
is slightly greater than 1 for complete, equal to 1 for complete multipartite,
and strictly less than 1 for other graphs. In the case of $\lambda_1 <1$
we gave an upper and lower estimate for the Cheeger constant of the graph
by $\lambda_1$  (see \cite{Bolla1}), illustrating that a smallest positive
normalized Laplacian eigenvalue
a separated from zero is an indication of the high edge-expansion
of the graph. 
In view of the above, this estimation is not valid for complete and
complete bipartite or multipartite graphs. Indeed, former ones are, in fact, 
super-expanders, while latter ones are so-called bipartite or multipartite 
expanders. By continuity, for large $n$, a $\lambda_1$ close to 1 
(from the left)  is also `suspicious', as it may indicate that our graph is
close to a bipartite or multipartite expander. The situation can even be more
complicated and also influenced by the upper end (near to 2 eigenvalues)
of the normalized Laplacian matrix, see \cite{Bol13}.
\end{itemize}

%Though, the open question of D. Stevanovic was answered, another open question
%emerges: whether the indefinite modularity and normalized modularity matrix
%of a graph, which is not complete multipartite, has the same number of negative
%and positive eigenvalues? Our conjecture is that yes, 
%but we leave the proof to the reader.

\end{document}